\documentclass{article}
\usepackage[utf8]{inputenc}
\usepackage[english]{babel}

\usepackage{biblatex}
\usepackage{hyperref}

\usepackage{ dsfont }
\usepackage{amssymb}
\usepackage{amsmath}        % rozšíření pro sazbu matematiky
\usepackage{amsfonts}       % matematické fonty
\usepackage{amsthm}         % sazba vět, definic apod.
\usepackage{bbding}         % balíček s nejrůznějšími symboly
\usepackage{bm}             % tučné symboly (příkaz \bm)
\usepackage{graphicx}       % vkládání obrázků
\usepackage{fancyvrb}       % vylepšené prostředí pro strojové písmo
\usepackage{indentfirst}    % zavede odsazení 1. odstavce kapitoly
\usepackage{icomma}         % inteligetní čárka v matematickém módu
\usepackage{dcolumn}        % lepší zarovnání sloupců v tabulkách
\usepackage{booktabs}       % lepší vodorovné linky v tabulkách
\usepackage{paralist}       % lepší enumerate a itemize
\usepackage{xcolor}         % barevná sazba
\usepackage{ dsfont }
\usepackage{textcomp}
\theoremstyle{plain}
\newtheorem{veta}{Věta}
\newtheorem{Thm}[veta]{Theorem}
\newtheorem{Prop}[veta]{Proposition}
\newtheorem{Example}[veta]{Example}
\newtheorem{Cor}[veta]{Corollary}

\theoremstyle{plain}
\newtheorem*{Def}{Definition}

\newtheorem{lemma}{Lemma}

\theoremstyle{remark}

\theoremstyle{plain}

\newenvironment{dukaz}{
  \par\smallskip\noindent
  \textit{Proof}.
}{
\newline
\rightline{$\qedsymbol$}
}

\usepackage{blindtext}
\title{RESTRICTED MINIMUM CONDITION IN REDUCED COMMUTATIVE RINGS}
\author{Dominik Krasula \thanks{This work is a part of  projects SVV-2020-260589 and PRIMUS/21/SCI/014}
}

\begin{document}
 \maketitle 
\begin{abstract} We say that a commutative ring $R$ satisfies the restricted minimum (RM) condition if for all essential ideals $I$ in $R$, factor  $R/I$ is an Artinian ring. We will focus on Noetherian reduced rings because in this setting known results for RM domains generalize well. However, as we will show, RM rings need not be Noetherian and may have nilpotent elements.

One of the classic results in the theory of RM rings is that for Noetherian domains RM condition corresponds to having Krull dimension at most one. We will show that this can be generalized to reduced Noetherian rings, thus proving that affine rings corresponding to curves are RM. We will give examples showing that the assumption that the ring is reduced is not superfluous.     

   We will prove that CDR domains are RM and this will allow us to give a new characterization of Dedekind domains.  Examples of RM rings for various classes of rings will be given. In particular, we will show that a ring of polynomials $R[x]$ is RM if and only if $R$ is reduced Artinian ring. And we will study the relation between RM rings and UFDs.

\end{abstract}
\textbf{Keywords:} restricted minimum condition, reduced commutative rings, noetherian rings, affine rings, polynomial rings, CDR  domains

\smallskip

\noindent \textbf{Mathematics Subject Classification:} 16P70, 	13E05; Secondary: 	13B25  , 	13F05

\smallskip

\noindent \textbf{Author:} Bc. Dominik Krasula

Charles University, Faculty of Mathematics and Physics, Prague, Czech Republic

resitelmam@gmail.com

ORCID: 0000-0002-1021-7364

\section{Introduction} For an associative unital ring $R$ we say that a right $R$-module $M$ is \textit{Artinian} (\textit{Noetherian}) if its lattice of submodules satisfies the descending (ascending) chain condition.  We say that $R$ is right Artinian (Noetherian) if it is Artinian (Noetherian) as a right $R$-module.

Unlike  general lattice, the descending chain condition on right ideals implies the ascending chain condition, as is shown by the classical Hopkins-Levitzki theorem \cite[Thm 15.20]{CatMod}. 
For the purposes of this article, it is enough to consider the commutative version due to Akizuki. It can be proved as a special case of Hopkins-Levitzky but it is historically older. Among other sources, direct proof can be found in \cite{Coh}. Proof using the more general Hopkins-Levitzky theorem can be found in appendices.
\begin{Thm}[Akizuki]\label{HP}
A commutative ring $R$ is Artinian iff it is right Noetherian and all prime ideals are maximal.
\end{Thm}
In particular, Artinian rings have composition series. In some settings, the Artinian condition is too strong. Reduced commutative rings are Artinian only in the trivial case, i.e., if they are a finite sum of fields  \hyperref[Akizuki]{[prop. \ref{Akizuki}]}. It is well known that a module is Artinian iff all of its factors are Artinian. The restricted minimum condition generalizes Artinian modules, requiring only  \textit{some} of the factors to be Artinian. 

This concept was first used during the 1950s in the work of I. S. Cohen\label{Cohen} \cite{Coh},
who studied commutative rings whose factors by \textit{non-zero} ideals are Artinian rings.  
Ornstein \cite{Ornstein} generalised this concept for non-commutative rings. In the 2000s, Király studied this condition for group rings \cite{Király}. 

 We will call domains with this condition \textit{RM domains}. This concept proved itself to be a useful tool in the study of Dedekind domains as it allowed both for proofs of new theorems and simpler proofs of known facts - see \cite{Coh} for details. We will use this method  to prove that a domain is CDR iff it is Dedekind \hyperref[Dedekind]{[Thm.~\ref{Dedekind}]}.

By Akizuki, all factors of RM domains are also Noetherian, hence RM domains are Noetherian too. Cohen proved  that domains of finite rank are RM and for local  domains, the converse is also true \cite[Thm. 10, 9]{Coh}. He also observed the following corollary of Akizuki's theorem:
\begin{Cor}\label{Cor} 
A Noetherian domain is RM iff its dimension is at most one, i.e., non-zero prime ideals are maximal 
\end{Cor}

However, in many aspects Cohen's definition is insufficient. If such a commutative ring  has zero-divisors, it has to be Artinian ring   \cite[Cor. 1]{Coh}. We also see that such rings are not closed under finite direct sums. 

In a module-theoretic setting, a module is said to satisfy the  restricted minimum condition if its factors by \textit{essential} submodules are Artinian modules. A~submodule is essential if it has a non-zero intersection with all non-zero submodules. This concept was studied, among others, by \cite{Huynh},  \cite{Breaz}, \cite{Kosan} and~\cite{Karami}.

\smallskip

All rings in this paper are considered to be unital commutative rings with $1\neq 0$.  From now on,  right Artinian rings will be called \textit{Artinian}. 

\begin{Def} 
We say that a ring satisfies  \emph{the restricted minimum condition}, or that it is an \emph{RM ring}, if all factors by essential ideals are Artinian rings, i.e., if it is an RM module as a module over itself. 
\end{Def}

RM rings are easily seen to be closed under taking factors and finite direct products. Any finitely generated module over an RM ring is RM. See \cite{Karami} for details. 

Important examples of RM rings (even in non-commutative setting) are Noetherian Hereditary rings \cite[Thm. 2.5]{Chatters}. However, RM ring needs not to be Noetherien, \hyperref[Counter]{[subsection 2.2]} and can have infinite global dimension \hyperref[-s]{[example \ref{-s}]}.

Any ideal containing a regular commutative element (i.e., element that is not a zero divisor) is an essential ideal. We see that for a domain our definition is equivalent to \hyperref[Cohen]{Cohen's}. Hence the study of RM domains is essentially a study of Noetherian domains of dimension 1. We will see that the corollary  \ref{Cor} generalizes to reduced rings \hyperref[reduced]{[Thm. \ref{reduced}]} hence affine rings corresponding to algebraic curves all possess the RM property \hyperref[affine]{[Cor. \ref{affine}]}.

Over the last decade, there was significant progress in the theory of modules over RM domains, generalizing some of the known results about torsion self-small abelian groups, see \cite{Breaz}, and \cite{Albrecht}. Some of these results were later generalized to commutative rings by Kosan and Žemlička \cite{Kosan}.
\begin{Thm}\cite[Thm 3.7]{Kosan}
The following are equivalent for a commutative ring $R$.

(i) $R$ is an RM ring.

(ii) If $M$ is singular $R$ module, then $M=\bigoplus_{P\in Max(R)} M_{[P]}$, where $M_{[P]}$ denotes sum of all finite-length submodules of $M$ whose composition factors are isomorphic to $R/P$.

(iii) $R/Soc(R)$ is Noetherian and every self-small singular module is finitely  generated. 
\end{Thm}

The aim of this paper is to study RM rings in terms of their internal properties, continuing the work of Cohen with the aid of results from the general theory of RM modules, mainly \cite{Kosan}.  We accompany the theory with constructions of several examples, focusing on polynomial rings and their factors.

\textbf{In section \ref{Reduced}}, we will generalize  \hyperref[Cohen]{Cohen's} characterization of RM domains to \hyperref[reduced]{reduced Noetherian rings}. \hyperref[polynomial]{We will  show that polynomial ring $R[x]$ is RM iff $R$ is reduced Artinian}. And we will \hyperref[Counter]{construct an example of a reduced RM ring that is not Noetherian.}

\textbf{In section \ref{Domain}}, we will study the relation of  RM domains  to various important classes of domains such as Dedekind domains, unique factorization domains (UFDs) and principal ideal domains (PIDs). \hyperref[CDR]{We will prove that containment-division  domains are RM domains.}  Then we will utilize this fact \hyperref[Dedekind]{in proving that Dedekind domains can be characterized as CDR domains.}  \hyperref[SectionPID]{We will show that GCD domains and domains of Laurent polynomials  are  RM domains iff they are principal ideal domains.}

\section{Reduced RM rings} \label{Reduced}
 For an ideal $I$, we will call the intersection of all prime ideals containing~$I$ \textit{the radical of I}, denoted by $\sqrt{I}$. We say that a ring is \textit{reduced} if it has no non-zero nilpotent elements, i.e., $\sqrt{0}=0$.  We say that a ring is a \textit{domain} if it has no non-zero zero divisors, i.e.,  the zero ideal is prime. 
 
 In this section, we study the RM condition for reduced RM rings.  We start by formulating some of the theory needed. In \hyperref[Reduced1]{subsection 2.1} we will present our results on Noetherian reduced rings. \hyperref[Counter]{In subsection 2.2} we construct an example of a reduced RM ring that is not Noetherian.

For most of this section, we will study rings using their lattice of ideals. For the convenience of the reader, we formulate some basic properties that we will implicitly use throughout the text.
\begin{lemma}\label{lattice}
Let $R$ be a ring and $I$ an ideal in $R$. Denote $\mathcal{L}(R)$  the lattice of ideals of $R$ and $[I, R]$ its sublattice of ideals containing $I$.

(1) There is a lattice  isomorphism  $\phi : [I,R] \rightarrow \mathcal{L}(R/I)$.

(2) Restricting $\phi$ on prime (maximal) ideals gives  a bijection between prime

~ ~ (maximal) ideals of $R$ containing $I$ and prime (maximal) ideals of $R/I$.

(3) If $E$ is an essential ideal in $R/I$ then $\phi^{-1}(E)$ is essential in $R$. 
\end{lemma}
Recall that factors by prime ideals are domains. As a corollary to $(3)$, all ideals containing a prime ideal are essential. The minimal primes themselves may or may not be essential:
\begin{Example} In domains, minimal prime (the zero ideal) is not essential.

Minimal primes contain nilradical $\sqrt{0}$. If it is a non-zero ideal, it is an  essential ideal, hence also minimal primes are essential. 

In $\mathds{Z}_5[x]/(x^2)$ the nilradical is the minimal prime ideal.  In $\mathds{Z}_6[x]/(x^2)$ the nilradical is not prime, consider $[x+2]_{(x^2)}\cdot [x+3]_{(x^2)}=[5x]_{(x^2)}$, so its minimal primes strictly contain an essential ideal.

\end{Example}

\phantomsection \label{dimension}Consider a prime ideal $P\leq R$, we define its \textit{height} to be the supremum of lengths of chains of prime ideals strictly included in $P$. We define \textit{Krull dimension of $R$}, abbreviated as $Kdim(R)$,  to be the supremum of heights of its primes. We define (Krull) dimension of  an ideal $I$ to be $Kdim(R/I)$.

Fields  have dimension 0. Dedekind domains that  are not fields have dimension one. If $R$ is a ring and $R[x]$ its ring of polynomials, then $Kdim(R[x])=1+Kdim(R)$ \cite[Cor. 10.13]{Eisenbud}. 

Krull dimension of RM rings is at most one.  Take any prime ideal $P$ in an RM ring $R$.  RM rings are closed under taking  factors,  so $R/P$ is  an  RM  domain  thus  all  non-zero  prime  ideals  are  maximal. Hence~the chain of primes containing $P$ can have length at most 1.

We say that a ring $R$ is \textit{semisimple} if $Soc(R)=R$. As we are restricted only to commutative rings, this is equivalent  to being a finite direct product of fields. We will now formulate special cases of the \hyperref[HP]{Akizuki  theorem}.
\begin{Prop}\label{Akizuki}
Let $R$ be a ring. 

(1) If $R$ is reduced then $R$ is Artinian iff it is a finite direct product of fields. 

(2)\label{(2)} If $R$ is a domain then $R$ is Artinian iff it is a field. 
\end{Prop}
\begin{dukaz}Domains have no non-trivial direct summands so $(2)$ follows from $(1)$. 

 Suppose $R$ is reduced and Artinian. The Jacobson radical $J(R)$ of an Artinian ring is known to be nilpotent.  Zero is the only nilpotent element in a reduced ring, hence $J(R)$ is the zero ideal.
 
 Because $R$ is of Krull dimension zero, its maximal ideals are also minimal. Noetherian rings have only finitely many minimal primes \cite[Thm. 88]{Kaplan}. Thus Artinian  rings are semilocal, i.e.,  there are only finitely many maximal ideals in $R$. We will denote them $M_1,\dots, M_n$. By the Chinese-Remainder Theorem $R/J(R)=R/M_1\times \dots \times R/M_n$. Hence  $R/J(R)\cong R$ is a finite product of fields. 
\end{dukaz}
\vfil \break 
\subsection{Reduced Noetherian Rings} \label{Reduced1}

As observed by \hyperref[Cohen]{Cohen}, RM domain is necessarily Noetherian. This is not true for rings with zero divisors, as shown by example \ref{Counter}. Noetherian RM domain is either field or of Krull dimension 1, and these are also sufficient conditions. We will now generalize this observation to reduced rings. The assumption that ring is reduced is not superfluous, as shown by the example~\ref{nilpotent}.

\begin{Thm}\label{reduced}
Suppose $R$ is a Noetherian reduced ring.  

Then $R$ is an RM ring iff $Kdim(R)\leq 1$.
\end{Thm}
\begin{dukaz}We have already seen that any RM ring has dimension at most one. To prove the other implication, we will show that minimal prime ideals  are not essential.  Thus any factor by an essential ideal is a Noetherian ring of dimension at most zero, hence it is Artinian by \hyperref[HP]{the Akizuki  theorem}

 Because $R$ is a reduced ring, the zero ideal is equal to the intersection of all  prime ideals.  Dimension  of  $R$ is finite, so it is enough to consider minimal primes in this intersection. Because $R$ is Noetherian, there are only finitely many minimal primes. (\cite[Thm. 88]{Kaplan}). We will denote them $P_1, \dots, P_n$.

 We have $0=\cap_{i\leq n} P_i$. We will show that for any $i\leq n$, the ideal $P_i$ is not essential.
 
 Set $Q_i:=\cap_{j\neq i}P_j$.  Because $Q_I\cap P_i = 0$, it is enough to show  that $Q_i\neq 0$.
 
 If $Q_i=0$ then $Q_i\subseteq P_i$. Since $P_i$ is a prime ideal, this implies that there is $j\in\{1,\dots, n\}$ such that $j\neq i$ and $P_j\subseteq P_i$. However, $P_i$ is minimal prime so  $P_i=P_j$, contradicting $i\neq j$.
 \end{dukaz}
\indent Recall that  $I\leq R$ is a radical ideal iff the factor  $R/I$ is a reduced ring. Using Hilbert's Basis theorem, we obtain the following consequence of the Theorem \ref{reduced}, showing a relation between RM rings and algebraic curves.  
\begin{Cor}\label{affine}
Let $R$ be a Noetherian ring and $R[\bar{x}]$ be a polynomial ring in finitely many variables over $R$. Take a radical ideal $I\leq R[\bar{x}]$  of dimension at most one.

Then the factor $R[\bar{x}]/I$ is an RM  ring. 
\end{Cor}
Take a ring $R$ and consider $R[x]$, it's ring of polynomials.  By  Hilbert Basis Theorem, $R[x]$ is Noetherian iff $R$ is. We now aim to characterize when $R[x]$ is an RM ring. 
\begin{Prop}\label{polynomial}Polynomial ring $R[x]$ is an RM ring iff $R$ is a reduced Artinian ring.
\end{Prop}
\begin{dukaz}Take a field $F$. Then $F[x]$ \hyperref[dimension]{has dimension 1}. $F[x]$ is Noetherian and reduced, hence by the proposition \ref{reduced} an RM ring. 

For any two rings $S$ and $T$, there is an isomorphism $S[x]\times T[x] \cong (S\times T)[x].$ Because RM rings are closed under finite direct products, we can conclude that if $R$ is a finite product of fields, then $R[x]$ is an RM ring. 

For the opposite direction, suppose that $R[x]$ is an RM ring. If  $R$ contains a proper essential ideal $E$, then $E[x]$, the ideal of polynomials with coeficients in $E$, is a proper essential ideal in $R[x]$. By the isomorphism $R[x]/E[x]\cong (R/E)[x]$, we conclude that $(R/E)[x]$ is an Artinian ring. However, polynomial rings are never Artinian as shown by  chain of ideals $(x)\supsetneq (x^2)\supsetneq (x^3)\supsetneq \dots \supsetneq 0$. 

We see that $R$ has no proper essential ideals. It follows that  $Soc(R)=R$. Otherwise, there would be a maximal ideal containing a socle. But such an ideal is clearly proper and  essential. We see that $R$ is semisimple, hence a finite direct product of fields.
\end{dukaz}
\begin{Example}\label{nilpotent}
Let $F$ be a field and denote $R:=F[x]/(x^2)$.

By the proposition \ref{polynomial}, $F[x]$ is an RM ring. Because $x^2$ is a regular element, we conclude that $R$ is an Artinian ring, hence  $Kdim(R)=0$. Polynomial ring $R[y]$ is then Noetherian of dimension 1. 

However, $R$ contains nilpotent elements, for example $[x]_{(x^2)}$, so by \hyperref[polynomial]{the previous proposition}, $R[y]$ is not an RM ring.
\end{Example}
\subsection{RM rings without the maximal condition}\label{Counter}
As discussed above, any RM domain is Noetherian. We will now show that this is not true for rings with zero divisors.  

Kosan and Žemlička showed that if $R$ is an RM ring then  $R/Soc(R)$ is Noetherian \cite[Thm. 3.4]{Kosan}. As a corollary, we see that an RM ring $R$ is Noetherian iff $Soc(R)$ is Noetherian as an $R$-module. Observe, that if  $Kdim(R)=0$ and   $R/Soc(R)$ is Noetherian then $R$ is RM by the \hyperref[HP]{Akizuki  theorem}.
\begin{Example}
Consider $R$ to be the set of all finite and cofinite subsets of $\mathds{N}$. 

Define addition on $R$ as  symmetric difference and define the product of two sets to be their intersection. Unary operator $-$ is defined as the identity, the empty set is 0 and the whole $\mathds{N}$ is 1.

Then $R$ has a structure of a reduced RM ring that is not Noetherian. 
\end{Example}
%\par\smallskip\noinden  \textit{Proof.$~$}
\begin{dukaz}
 For any set $S\in R$, we can define ideal $I(S)$ to be the set of all finite and cofinite subsets of $S$. We see that $R$ is not Noetherian, as we can consider the infinite ascending chain $I(\{1\})\subsetneq I(\{1,2\})\subsetneq I(\{1,2,3\})\subsetneq \dots \subsetneq R$.

Observe that for any ideal $J\leq R$, if $S\in J$ then $I(S)$ is an ideal contained in $J$. Hence  simple ideals of $R$ are the ideals generated by singletons. We conclude that $Soc(R)$ is the  ideal consisting of all finite sets. 

The RM condition in $R$ follows from $R/Soc(R)$ being the two-element ring. Let $\pi$ be the canonical projection $R\mapsto R/Soc(R)$. If $F\in R$ is a finite set, then it is contained in $Soc(R)$, so $\pi(F)=0$. Now take two cofinite sets $X,Y\in R$. Then $X+Y$ is a finite set, hence $X-Y=X+Y\in Soc(R)$ so $\pi(X)=\pi(Y)$.
\end{dukaz} %\newline\rightline{$\qedsymbol$}
\section{RM domains}\label{Domain}
A domain is Artinian iff it is a field. The concept of an RM domain gives a weaker condition. Infinite strictly descending chains of ideals are allowed, but only if the intersection of all ideals in the chain is the zero ideal. 
%Equivalently, domain $R$ is an \textit{RM domain} if for a non-zero ideal $I\leq R$ the factor $R/I$ is an Artinian ring. 
\begin{Example}\label{PID}
Let $R$ be a  principal ideal domain (PID).

Take $a\in R$ and an infinite descending chain $(a_1)\supseteq (a_2)\supseteq \dots \supseteq (a)$. \newline If $a\neq 0$ then this chain stabilizes.  Take $r_i\in R$ such that $a_ir_i=a$. $R$ is PID, in particular is is Noetherian. Hence the ascending chain $(r_1)\subseteq (r_2)\subseteq \dots \subseteq R$ stabilizes, so the origininal chain stabilizes too. We see that PIDs are RM. 
\end{Example}
Later, we will generalize this observation to CDR domains.  

We say that a ring is \textit{semilocal} if it has  finitely many maximal ideals. \hyperref[(2)]{We have seen that Artinian rings are semilocal}. Let $D$ be a domain that is not a field. By the Chinese-Remainder Theorem, if $D$ is semilocal then $J(D)\neq 0$. The converse is not true in general but holds for RM domains.  If $J(D)\neq 0$ in an RM-domain then the Artinian factor $D/J(D)$ is semilocal, hence $D$ itself is semilocal. 
\subsection{CDR domains}\label{SectionCDR}
We say that a domain $R$ is a \textit{Dedekind domain} if every proper ideal can be written as a product of finitely many (not necessarily distinct) maximal ideals. Take $I\leq R$ an ideal. Because $R$ is RM, there is only finitely many maximal ideals containing $I$. It follows that factors of Dedekind domains have only finitely many ideals.

RM domains appeared  as a generalization of Dedekind domains in several  contexts.  Cohen \cite{Coh} characterizes Dedekind domains as integrally closed RM domains. Recalling that Dedekind domains can be characterized as (Noetherian) hereditary domains, the RM condition for Dedekind domains follows from the work of Chatters \cite[Thm. 2.2]{Chatters}. Albrecht and Breaz characterized torsion modules over RM domains, generalizing characterization of torsion modules over Dedekind domains \cite[Thm. 6]{Breaz}. 

Let $I$, $J$ be two ideals. We say that $I$ \textit{divides} $J$, denoted $I \mid J$, if there is  ideal $H$ such that $IH=J$. A ring will be called \textit{containment division ring}, or \textit{CDR} for short, if for any proper ideals $I\supseteq J$, it holds that $I\mid J$. Sometimes in literature, CDRs are called \textit{multiplication rings.}  

It is well-known that Dedekind domains are CDR. In the theorem \ref{Dedekind}, we will show that this property actually characterizes Dedekind domains. This proposition sometimes appears in the literature, see \cite[Thm. 3.10.]{Multiplicative}, and can be seen as a consequence of \cite[Thm. 6]{Breaz}. We however, aim to give a simpler proof using the notion of an RM domain.

For a domain $R$, denote $F$ its field of fractions. An $R$-submodule $M$ of $F$ will be called \textit{fractional ideal} if there exists an element $r\in R$ such that $rM\subseteq R$. Ideals of $R$ are fractional ideals. We say that a fractional ideal $M$ is invertible if there exists another fractional ideal $N$ such that $MN=R$. 
\begin{Prop}\label{CDR}
CDR domains are RM domains
\end{Prop}
\begin{dukaz}Consider a CDR domain $R$ and $F$ its field of fractions. For any non-zero ideal $I$ in $R$ there exists a principal ideal $0\neq (i)$ such that $(i)\subseteq I$. Because $R$ is a CDR, there exists the ideal $H$ such that $HI=(i)$. The element $i$ is regular, hence also invertible in $F$. This means that the ideal $(i)$ is invertible. Since multiplication of ideals is associative, we see that $I$ is invertible too. 

We have seen that all non-zero ideals of $R$ are invertible as a fractional ideals. It follows that any  ideal $I\leq R$ is finitely generated. Consider $J$ fractional ideal, the inverse of I. There exist $a_i\in I$ and $b_i\in J$ 
such that $1=a_1b_1+\dots +a_nb_n$. If we multiply the equation by $a\in I$ we get $a=a_1(ab_1)+\dots+ a_n(ab_n)$. Because $I$ is the inverse ideal to $J$, elements $ab_i$ lie in $R$. We obtain that any $CDR$ domain is Noetherian. 

\phantomsection \label{NoeCDR}
Now consider a decreasing chain of ideals $R\supseteq I_0\supseteq I_1 \supseteq I_2 \dots \supseteq I\neq 0 $

$I_{i+1}\subseteq I_i$ so by $R$ being a CDR there exists the ideal $H_i$ such that $I_{i+1}=I_iH$. Denote $J_i$ the ideal of $R$ such that $J_iI_i=I$. We have 
\[I_iJ_i=I_{i+1}J_{i+1}=I_iHJ_{i+1}\]

Ideal $I_i$ is invertible so we obtain $J_i=HJ_{i+1}$ hence 
$J_0\subseteq J_1\subseteq J_2\dots\subseteq R$. Because $R$ is  Noetherian, this chain stabilizes so $I_0\supseteq I_1 \dots \supseteq I$ stabilizes too. 
\end{dukaz}
\begin{Thm}\label{Dedekind}
A domain $R$ is a Dedekind domain iff it is a CDR.
\end{Thm}
\begin{dukaz}Suppose $R$ is a CDR domain and consider a non-zero ideal $I$ in $R$.  We want to write $I$ as a product of finitely many maximal ideals. If $I$ is maximal we are done. If $I$ is  not maximal, then there exists a maximal ideal $P_1$ containing $I$. Since $R$ is a CDR there is $I_1$ such that $P_1I_1=I$. If $I_1$ is not maximal we can find a maximal ideal $P_2$ and ideal $I_2$ such that $P_1P_2I_2=I$ and so on. 

This process necessarily stops, giving us the desired decomposition. Suppose this is not true. Then we have an infinite descending chain $P_1\supseteq P_1P_2\supseteq P_1P_2P_3 \dots \supseteq I$. Since any CDR domain is an RM domain we see that the factor $R/I$ is Artinian so the chain stabilizes.

On the other hand, suppose that $R$ is a Dedekind domain and we have two non-zero ideals $I\subseteq J$. Both ideals have decomposition into maximal ideals $I=P_1\dots P_n\subseteq Q_1 \dots Q_m=J\subseteq Q_1$. The ideal $Q_1$ is prime hence there is $i$ such that $P_i$  is contained in $Q_1$ and by maximality of $P_i$ we have $P_i=Q_1$. We can assume that $i=1$ and by \cite[Cor. 11.9.]{Eisenbud}  ideal $Q_1$ is invertible. 
We get $P_2\dots P_n\subset Q_2\dots Q_m$. Repeating the argument we get that $J=I\cdot Q_{n+1}\dots Q_m$, hence $R$ is a CDR. 
\end{dukaz}
 \subsection{Principal Ideal Domains}\label{SectionPID}
In this section, we turn our attention to unique factorization domains (UFDs) and PIDs. By the proposition \ref{CDR}, PIDs are RM domains.  Because a UFD needs not to be Noetherian  (take a ring of polynomials in infinitely many variables over UFD) it is clear that there exist UFDs that are not RM domains.  \hyperref[UFD]{We will show examples}  of Noetherian UFDs that are not RM domains and examples of RM domains that are not UFDs.

 We say that a domain is a GCD if for any two regular elements the greatest common divisor exists.  Recall that a Dedekind domain is a GCD iff it is a PID.  This can be easily generalized to  RM domains. If $R$ is GCD and RM it is, in particular, Noetherian GCD  thus it is UFD. All UFDs are integrally closed \cite[Thm 2.3]{Cohn} so we have proved:
  
\begin{Prop}\label{GCD}
Let R be a GCD. Then it is a  PID iff it is an RM domain.
\end{Prop}

One can prove the above proposition without dealing with integrall close-ness. This more direct proof can be found in appendices.

\indent Király \cite[Thm. 1.2]{Király} showed that a ring of Laurent  polynomials over a field is an RM domain. We now prove the converse of this observation.  
\begin{Prop}
Let $R$ be a domain. The ring of Laurent polynomials $R[x,x^{-1}]$ is an RM domain iff it is a PID which is equivalent with $R$ being a field.
\end{Prop}
\begin{dukaz}If $R$ is a field, polynomial ring $R[x]$ is a PID. Ring $R[x,x^{-1}]$ can be viewed as the localization of $R[x]$ at $x$ so it is also a PID, hence an RM domain. 

Suppose that $R[x,x^{-1}]$ is an RM domain. Denote $I$ the Kernel of the evalutation homomorphism $f: R[x,x^{-1}] \rightarrow R$ sending $x$ to the identity element in~$R$. Then $R[x,x^{-1}]/I\cong R$. We see that $R$ is an Artinian domain hence a field. 
\end{dukaz}\indent We will now utilize the Corollary \ref{affine}, to get some concrete examples of RM domains. Recall that $R[x]$ is a UFD if and only if $R$ is a UFD. 
\begin{Cor}
Let $R$ be a Dedekind domain and $f\in R[x]$ an irreducible monic polynomial. Then $R[x]/(f)$ is an RM ring.  
\end{Cor}
\begin{dukaz}By the theorem \ref{Dedekind}, Dedekind domain is an RM domain hence  $Kdim(R)\leq 1$. Dimension of the polynomial ring  $R[x]$ is then at most 2.  $R[x]$ is a domain so non-zero prime ideals in $R[x]$ have dimension at most 1.

Dedekind domains are integrally closed. So any principal ideals generated by monic irreducible polynomials are prime ideals \cite[Cor. 4.12]{Eisenbud}.  By the Corollary \ref{affine}, factors by such prime ideals are RM rings.
\end{dukaz}
\begin{Example}\label{-s}
Ring of integers $\mathds{Z}$ is a PID, hence a Dedekind domain.

Consider $s\in \mathds{N}$ that is not a square of an integer. 

Then ring $\mathds{Z}[\sqrt{-s}]\cong \mathds{Z}[x]/(x^2+s)$ is an RM domain. 

It is known that some of these rings are not integrally closed, e.g. $\mathds{Z}[\sqrt{-3}]$ so we obtain examples of RM domains that are not Dedekind.
\end{Example}
\begin{Example}\label{UFD}~

(1) Ring $\mathds{Z}[x]$ is a Noetherian UFD. By Prop. \ref{polynomial} it is not an RM domain.

(2) By example \ref{-s}, ring $\mathds{Z}[\sqrt{-14}]$ is an RM domain. 

It is not a UFD as $15=3\cdot 5 =(1+\sqrt{-14})(1-\sqrt{-14})$.

(3) Ring $\mathds{Z}[\sqrt{-14}][x]$ is Noetherian. It is neither a UFD nor an RM domain. 
\end{Example}

\smallskip
\noindent\textbf{Aknowledgements:}
 Author wishes to express his thanks and appreciation to Doc. Jan Žemlička for introducing him to the concept of an RM ring and for his guidance and valuable suggestions during the research. 
 
 \vfil \break

  \section{Appendices}
 
 \subsection{Akizuki theorem as a special case of Hopkins-Levitzki}
 
 \begin{Thm}[Hopkins-Levitzki]
A ring $R$ is right Artinian iff it is right Noetherian, the Jacobson radical $J(R)$ is nilpotent and $R/J(R)$ is a semisimple ring. 
\end{Thm}

\textbf{Suppose that R is Artinian} Let  $P\leq R $ be a prime ideal. Then $R/P$ is an Artinian domain, thus a field. In particular,  $P$ is a maximal ideal. Thus we arrive at $Kdim(P)=0$. The fact that $R$ is Noetherian follows directly from the Hopkins-Levitzki.

\textbf{Suppose that $R$ is Noetherian of dimension 0}. Then the only prime ideals in $R$ are maximal ideals, hence $J(R)$ equals the radical of the zero ideal. Because $R$ is Noetherian, $J(R)$ is finitely generated, so it is a nilpotent ideal.

We see that maximal ideals are minimal primes in $R$. Noetherian rings have only finitely many minimal primes \cite[Thm. 88]{Kaplan}. Hence there are only finitely many maximal ideals in $R$. We will denote them $M_1,\dots, M_n$. By the Chinese-Remainder Theorem $R/J(R)=R/M_1\times \dots \times R/M_n$. Hence  $R/J(R)$ is a semisimple ring and $J(R)$ is nilpotent. We conclude that $R$ is Artinian.
\newline 
\rightline{$\qedsymbol$}
\subsection{Direct proof of the proposition \ref{GCD}}
Suppose $R$ is a GCD domain. We will show that it is RM if and only if it is a PID.

\begin{dukaz}By the example \ref{PID}, all PIDs are RM so this  implication is clear. 

Suppose $R$ is both a GCD and an RM domain. RM domains are Noetherian so in particular, $R$ satisfies the ascending chain condition on principal ideals so $R$ is a UFD. 

Prime ideals are prinicipal in $R$. To prove it, take $P\leq R$ a prime ideal. Take non-zero  $p\in P$ and consider its factorization into prime elements $p=p_1^{k_1}\cdot ...\cdot p_n^{k_n}$. By $P$ being prime, there exists $i$ such that $(p_i)\subseteq P$. Recall that for an RM domain, prime ideals are maximal hence $P=(p_i)$. 

By Zorn's lemma, any proper ideal is contained in some maximal ideal. We have shown that maximal ideals are principal. So all elements in proper ideal have a common divisor. Thus any set of coprime elements generates the whole~$R$.

Now take a general non-zero proper ideal $I$. Mark $g_1,\dots, g_n$ its  generators and let $g$ be their greatest common divisor. Consider set $\{r_1, r_2\dots r_n\}\subseteq R$ such that $gr_i=g_i$. We see that $(g)(r_1,\dots,r_n)=I$. Elements  $r_1,\dots,r_n$ are coprime so they generate the whole $R$ thus $(g)=I$.
\end{dukaz}

\vfil \break

\printbibliography

\end{document}